\input amstex
\documentstyle{amsppt}

\topmatter
\centerline {\bf AN ARITHMETIC FORMULA FOR CERTAIN COEFFICIENTS OF}
\title
the Euler product of Hecke polynomials
\endtitle
\author
Xian-Jin Li
\endauthor
\abstract
 In 1997 the author [11] found a criterion for the Riemann 
hypothesis for the Riemann zeta function, involving
the nonnegativity of certain coefficients associated with
the Riemann zeta function.  In 1999 Bombieri and Lagarias [2]
obtained an arithmetic formula for these coefficients
using the ``explicit formula'' of prime number theory.
In this paper, the author obtains an arithmetic formula
for corresponding coefficients associated with the
Euler product of Hecke polynomials, which is essentially
a product of $L$-functions attached to weight $2$ cusp
forms (both newforms and oldforms) over Hecke congruence
subgroups $\Gamma_0(N)$.  The nonnegativity of these 
coefficients gives a criterion for the Riemann hypothesis
for all these $L$-functions at once. 
\endabstract
\address
Department of Mathematics, Brigham Young University,
Provo, Utah 84602
\endaddress
\email
xianjin\@math.byu.edu
\endemail
\thanks
Research supported by National Security Agency MDA904-03-1-0025
\endthanks
\endtopmatter
\document

\heading
1.  Introduction
\endheading

     We denote by $S_2(N)$ the space of all cusp forms
of weight $2$ with the principal character modulo $N$
for the Hecke congruence subgroup
$\Gamma_0(N)$ of level $N$.  That is,
$f$ belongs to $S_2(N)$ if, and only if,
$f$ is holomorphic in the upper half-plane, satisfies
$$f\left({az+b\over cz+d}\right)=(cz+d)^2 f(z)$$
for all $\left(\smallmatrix a&b\\ c&d\endsmallmatrix\right)
\in \Gamma_0(N)$, satisfies the usual regularity conditions at
the cusps of $\Gamma_0(N)$, and vanishes at each cusp
of $\Gamma_0(N)$.

   The Hecke operators $T(n)$, $n=1,2,\cdots$,
are defined by
$$(T(n) f)(z)={1\over n}\sum_{ad=n, (a,N)=1}a^2
\sum_{0\leq b<d} f\left({az+b\over d}\right) \tag 1.1$$
for any function $f\in S_2(N)$.
The Fricke involution $W$ is
defined by
$$(Wf)(z)=N^{-1}z^{-2}f(-1/Nz),$$
and the complex conjugation operator $K$ is defined by
$$(Kf)(z)=\bar f(-\bar z).$$
Set $\bar W=KW$.   A function $f\neq 0$ in
$S_2(N)$ is called a newform if
it is an eigenfunction of $\bar W$ and of all the Hecke
operators $T(n)$.

   Let $f$ be a newform in $S_2(N)$ normalized so
that its first Fourier coefficient is $1$.  Then it
have the Fourier expansion
$$f(z)=\sum_{n=1}^\infty \lambda(n) e^{2\pi i nz}$$
with the Fourier coefficients equal to the eigenvalues of
Hecke operators.  The Hecke $L$-function associated with
$f$ is given by
$$L_f(s)=\sum_{n=1}^\infty {\lambda(n)\over n^s}$$
for $\Re s>{3\over 2}$; see \S7.2 in Iwaniec [8].
It has the Euler product
$$L_f(s)=\prod_{p|N}(1-\lambda(p)p^{-s})^{-1}
\prod_{p\nmid N} (1-\lambda(p)p^{-s}+p^{1-2s})^{-1}.
\tag 1.2 $$
If we denote
$$\xi_f(s)=N^{s/2}(2\pi)^{-s}
\Gamma({1\over 2}+s)L_f({1\over 2}+s), \tag 1.3$$
then $\xi_f(s)$ is an entire function and satisfies the
functional identity
$$\xi_f(s)=w\xi_f(1-s)$$
where $w=(-1)^r$ with $r$ being the vanishing order of
$\xi_f(s)$ at $s=1/2$.

   Let $I$ be the identity map acting on the space $S_2(N)$.
For each prime $p\nmid N$, we put
$$H_p(u)=\det(I-T(p)u+pu^2I)$$
which we call a Hecke polynomial; see Ihara [7].  
The Euler product $L_N(s)$ of Hecke polynomials is defined by
$$L_N(s)=\prod_{p\nmid N}H_p(p^{-s})^{-1}, \tag 1.4$$
which is the main object of study in this paper.
Let
$$\xi_N(s)=N^{gs/2}(2\pi)^{-gs}\Gamma^g({1\over 2}+s)
L_N({1\over 2}+s) \tag 1.5$$
where $g$ denotes the dimension of the space $S_2(N)$.
It will be shown (Lemma 4.2) that $\xi_N(s)$ is a
product of the $\xi_f(s)$'s over a Hecke eigenbasis
of the space $S_2(N)$ times a finite Euler factor,
that $\xi_N(s)$ is an entire function, and
that zeros of $\xi_N(s)$ in the strip $0<\Re s<1$ appear in
pairs $\rho$ and $1-\rho$.

   Let
$$\tau_N(n)=\sum_\rho [1-(1-{1\over\rho})^{-n}] \tag 1.6 $$
for $n=1,2,\cdots$, where the sum on $\rho$ runs over all
zeros of $\xi_N(s)$ taken in the order given by
$|\Im \rho|<T$ for $T\to\infty$ with a zero of multiplicity
$\ell$ appearing $\ell$ times in the list.  If $\rho=0$ is
a zero of $\xi_N(s)$, then $(1-1/\rho)^{-n}$ in (1.6) is
interpreted to be $0$.

  In order to state the main result of this paper (Theorem 1.2), 
we need an explicit formula for the trace tr$T(p^k)$ of Hecke
operators $T(p^k)$ acting on the space $S_2(N)$
for all primes $p\nmid N$ and for $k=1,2,\cdots$.
This formula is given by the Eichler-Selberg
trace formula obtained in Oesterl\'e [16].
We denote by $\varphi$ the Euler
$\varphi$-function.  Let
$$\psi(N)=N\prod_{p|N}  (1+1/p).$$

\proclaim{Lemma 1.1}  Let $p$ be a prime with $(p, N)=1$.
Then we have
$$\aligned \text{tr}(T(p^k))&={1+(-1)^k\over 24}\psi(N)
  +{p^{k+1}-1\over p-1}\\
&-\sum_{t\in\Bbb Z, t^2<4p^k}\sum_{\underset{{t^2-4p^k\over m^2}
  \equiv 0, 1 (\mod 4)}\to {m\in \Bbb Z^+, m^2|t^2-4p^k}}{h((t^2-4p^k)/m^2)
  \over w((t^2-4p^k)/m^2)}{\psi(N)\over \psi(N/(N,m))}\mu_{t,p^k,m}\\
  &-{1\over 2}\sum_{0<d|p^k}\min(d, p^k/d)\sum_{(c, {N\over c})|
  (N, {p^k\over d}-d), c|N}\varphi((c, N/c))
   \endaligned $$
 for $k=1,2,\cdots$,
  where $\mu_{t,p^k,m}$ is the number of solutions $x$
  modulo $N$ of the equation $x^2-tx+p^k\equiv 0 (\mod N(N, m))$,
  and where $h(f)$ and $w(f)$ are the class number and
   the number of units in the ring of integers of the
  imaginary quadratic field of discriminant $f<0$ respectively.
\endproclaim

  Let $\gamma=0.5772\cdots$ be Euler's constant, let
$$\Lambda(k)=\cases \ln p, &\text{if $k$ is a
positive power of a prime $p$;}\\
0, &\text{otherwise,} \endcases $$
and let $d(k)$ be the number of positive divisors
of $k$.  For $m=1, 2, \cdots$, we denote by $\nu_m$
the dimension of the subspace generated by all newforms
in $S_2(m)$.  

   In this paper we obtain the following
arithmetic formula for the $\tau_N(n)$'s, which generalizes
an arithmetic formula of Bombieri and Lagarias [2] for the
Riemann zeta function.

 \proclaim{Theorem 1.2}  Let $\tau_N(n)$ be given in
 (1.6).  Then we have
$$\aligned \tau_N(n)=& {n\over 2}\ln(N^{\nu_N}
\prod_{\underset{1<m<N}\to{m|N}} m^{\nu_md(N/m)})\\
&-\sum_{l=1}^n \binom nl {(-1)^{l-1}\over (l-1)!}
\sum_{\underset{(m, N)=1}\to {m=1}}^\infty
{\Lambda(m)\over m^{3/2}}B(m)(\ln m)^{l-1}\\
&-ng\left(\ln(8\pi)+\gamma-2\right)
+g\sum_{m=2}^n\binom nm \sum_{l=1}^\infty
{(-1)^m\over (l+1/2)^m} \endaligned  $$
for all positive integers $n$, where
$B(p^k)= \text{tr}(T(p^k))-p\,\text{tr}(T(p^{k-2}))$
for $p\nmid N$ and where the second term on the right side of
the identity is interpreted as the limit
$$\lim_{X\to\infty} \sum_{l=1}^n \binom nl {(-1)^{l-1}\over (l-1)!}
\sum_{\underset{(m, N)=1}\to {m<X}}
{\Lambda(m)\over m^{3/2}}B(m)(\ln m)^{l-1}$$
which exists.
\endproclaim

In [2], Bombieri and Lagarias generalized a criterion of
the author for the Riemann hypothesis [11] and obtained
the following useful theorem.

\proclaim{Theorem 1.3} (Bombieri-Lagarias [2])
 Let $\Cal R$ be a set of complex
numbers $\rho$, whose elements have positive integral multiplicities
assigned to them, such that $1\not\in\Cal R$ and
$$\sum_\rho {1+|\Re \rho|\over (1+|\rho|)^2}<\infty.$$
Then the following conditions are equivalent:
\roster
\item $\Re\rho\leq {1\over 2}$ for every $\rho$ in $\Cal R$;
\item $\sum_\rho \Re[1-(1-{1\over\rho})^{-n}]\geq 0$
   for $n=1,2,\cdots$.
\endroster
\endproclaim

   As a corollary of Theorem 1.3 we obtain a criterion
for the location of all nontrivial zeros of Hecke $L$-functions
associated with all cusp forms which form an orthonormal basis
in $S_2(N)$ and which consist of eigenfunctions of all the
Hecke operators $T(n)$ with $(n, N)=1$.

\proclaim{Corollary 1.4}  All zeros of $\xi_N(s)$ in the
 strip $0<\Re s<1$ lie on the critical line
$\Re s=1/2$ if, and only if, $\tau_N(n)\geq 0$
for all positive integers $n$.
\endproclaim

   In a subsequent paper [13], the author generalized results of 
this paper to the space $S_k(N, \chi)$ of cusp forms of
weight $k$ for all integers $k>2$ and of character $\chi$
for all Dirichlet characters of modulus $N$ with
$\chi(-1)=(-1)^k$.
 
   This paper is organized as follows:  The Eichler-Selberg trace
formula for Hecke congruence subgroups, which is needed for
Theorem 1.2, is given in section 2.  An arithmetic formula
is obtained in section 3 for a sequence of numbers, whose
positivity implies that all zeros of Hecke $L$-functions
associated with newforms lie on the critical line.
 This formula will be used in the proof of
 Theorem 1.2.   In section 4 we give some preliminary
 results for the proof of Theorem 1.2.  Finally,
 Theorem 1.2 is proved in section 5.

   This research started while the author attended the Workshop on
Zeta-Functions and Associated Riemann Hypotheses, New York
University, Manhattan, May 29 - June 1, 2002.  He wants to thank
the American Institute of Mathematics, Brian Conrey, and Peter
Sarnak for the support.  He wishes to thank William Duke for 
his encouragement.  The author also wishes to thank the
referee for his valuable suggestions of improving the
presentation of this paper.

\heading
2.  The Eichler-Selberg trace formula
\endheading

  In this section, we state J. Oesterl\'e's formula [16] for the trace 
of Hecke operators acting on the space $S_2(N)$.

      Let $\chi_0$ be the principal character modulo $N$, and let
$\chi_0(\sqrt n)=0$ if $n$ is not the square of an integer.

  \proclaim{Lemma 2.1} (Th\'eor\`eme $3^\prime$, [16]) (cf. Cohen [3]) 
For every
  positive integer $n$, the trace tr$(T(n))$ of the Hecke operator $T(n)$
  acting on the space $S_2(N)$ is given by
  $$\aligned \text{tr}(T(n))&={1\over 12}\chi_0(\sqrt n)\psi(N)
  +\sum_{\underset{\gcd(N, n/d)=1}\to {0<d|n}}d\\
  &-\sum_{t\in\Bbb Z, t^2<4n}\sum_{\underset{{t^2-4n\over m^2}
  \equiv 0, 1 (\mod 4)}\to {m\in \Bbb Z^+, m^2|t^2-4n}}{h((t^2-4n)/m^2)
  \over w((t^2-4n)/m^2)}\mu(t, n, m)\\
  &-{1\over 2}\sum_{0<d|n}\min(d, n/d)\sum_{\underset{\gcd(c, {N\over c})|
  \gcd(N, {n\over d}-d)}\to {c|N}}\varphi(\gcd(c, N/c))\chi_0(y)
   \endaligned $$
  where the integer $y$ is defined modulo $N/\gcd(c, N/c)$ by
  $y\equiv d (\mod c), y\equiv n/d (\mod N/c)$, where
  $$\mu(t, n, m)={\psi(N)\over \psi(N/\gcd(N,m))}
  \sum_{\underset{x^2-tx+n\equiv 0 (\mod N\gcd(N, m))}\to
  {x (\mod N)}}\chi_0(x),$$
  and where $h(f)$ and $w(f)$ are respectively the class number
  and the number of units in the ring of integers of the
  imaginary quadratic field of discriminant $f<0$.
   \endproclaim

  \proclaim{Lemma 2.2}   Let $n=p^k$ for a prime $p\nmid N$ and
for a positive integer $k$, let $d|n$, and let $c|N$.  If $y$ is defined modulo
$N/(c, N/c)$ by $y\equiv d (\mod c)$ and $y\equiv n/d (\mod N/c)$,
then it is unique and $(y, N)=1$.\endproclaim

\demo{Proof}  Suppose that $q|(y, N)$ for a prime $q$.  Then
$q|c$ or $q|N/c$.  Assume that $q|c$.  Since $y\equiv d (\mod c)$,
we have $q|d$.  This contradicts to that $d$ is a power of $p$
with $(p, N)=1$.  The case when $q$ divides $N/c$ can be treated
similarly.  Thus, we have proved that $(y, N)=1$.

  Suppose that $y_1$ is another integer satisfying the hypotheses.
Then $y-y_1\equiv 0 (\mod c)$ and $y-y_1\equiv 0 (\mod N/c)$.
Let $\ell=(c, N/c)$.  We write $c=\ell c^\prime$,
$N/c=\ell N^\prime$ for some integers $c^\prime, N^\prime$.
Then $y-y_1\equiv 0 (\mod \ell c^\prime N^\prime)$.  That is,
$$y\equiv y_1 (\mod N/\ell).$$
Hence, $y$ is unique modulo $N/\gcd(c, N/c)$. \qed\enddemo

  \demo{Proof of Lemma 1.1}  Let $n=p^k$ for a prime $p\nmid N$ and
for a positive integer $k$.    Then we have
$${1\over 12}\chi_0(\sqrt n)\psi(N)={1+(-1)^k\over 24}\psi(N) \tag 2.1 $$
and
$$\sum_{d>0, (N, n/d)=1, d|n}d={p^{k+1}-1\over p-1}. \tag 2.2 $$

 Let $x$ be an integer given by the equation
$x^2-tx+n\equiv 0 (\mod N\gcd(N, m))$, where $m$ is an integer.
Since $n=p^k$ with $(p, N)=1$,  we have $(x, N)=1$.  Hence, we have
$$\chi_0(x)=1. \tag 2.3 $$
By (2.1), (2.2), (2.3), Lemma 2.1, and Lemma 2.2 we have
$$\aligned \text{tr}(T(p^k))&={1+(-1)^k\over 24}\psi(N)
  +{p^{k+1}-1\over p-1}\\
&-\sum_{t\in\Bbb Z, t^2<4p^k}\sum_{\underset{{t^2-4p^k\over m^2}
  \equiv 0, 1 (\mod 4)}\to {m\in \Bbb Z^+, m^2|t^2-4p^k}}{h((t^2-4p^k)/m^2)
  \over w((t^2-4p^k)/m^2)}{\psi(N)\over \psi(N/(N,m))}\mu_{t,p^k,m}\\
  &-{1\over 2}\sum_{0<d|p^k}\min(d, p^k/d)\sum_{(c, {N\over c})|
  (N, {p^k\over d}-d), c|N}\varphi((c, N/c))
   \endaligned $$
  where $\mu_{t,p^k,m}$ is the number of solutions $x$
  modulo $N$ of the equation $x^2-tx+p^k\equiv 0 (\mod N(N, m))$,
  and where $h(f)$ and $w(f)$ are the class number and
   the number of units in the ring of integers of the
  imaginary quadratic field of discriminant $f<0$ respectively.

  This completes the proof of the lemma. \qed\enddemo

\heading
3.  An arithmetic formula
\endheading

   In this section, an arithmetic formula is given
for the $\tau_f(n)$ (see (3.1) below) which generalizes
that of Bombieri and Lagarias [2] for the
Riemann zeta function.

   Let $f$ be a normalized newform in $S_2(N)$, and let
$\xi_f(s)$ be given in (1.3).  Put
$$\tau_f(n)=\sum_\rho [1-(1-{1\over\rho})^n] \tag 3.1 $$
for $n=1,2,\cdots$, where the sum is over all the zeros of $\xi_f(s)$
in the order given by $|\Im\rho|<T$ for $T\to\infty$
with a zero of multiplicity
$\ell$ appearing $\ell$ times in the list.

   Assume that $f$ is a normalized newform in $S_2(N)$.
For each prime number $p$, let $\alpha_p$ and $\beta_p$
be the two roots of $T^2-\lambda(p)T+p$ where $\lambda(p)$
is given in (1.2).  Put
$$b_f(p^m)=\cases\lambda(p)^m, &\text{ if $p|N$;}\\
\alpha_p^m+ \beta_p^m, &\text{ if $(p,N)=1$.}\endcases
\tag 3.2$$

  The following theorem is essentially obtained in Li [12],
where the theorem is presented in the context 
of $L$-series of elliptic curves.   For the convenience of
readers, we give a proof of the following theorem here which 
is a minor modification of that given in [12].

  \proclaim{Theorem 3.1}  Assume that $f$ is a normalized newform
in $S_2(N)$.  If $\tau_f(n)$ is given in (3.1), then we have
$$ \aligned \tau_f(n)=& n\left(\ln {\sqrt N\over 2\pi}-\gamma\right)
-\sum_{j=1}^n \binom nj {(-1)^{j-1}\over (j-1)!}\sum_{m=1}^\infty
{\Lambda(m)\over m^{3/2}} b_f(m)(\ln m)^{j-1}\\
&+n\left(-{2\over 3}+\sum_{l=1}^\infty {3\over l(2l+3)}\right)
+\sum_{j=2}^n\binom nj (-1)^j\sum_{l=1}^\infty
{1\over (l+1/2)^j} \endaligned $$
for $n=1,2,\cdots$, where the second term on the right side of
the identity is interpreted as the limit
$$\lim_{X\to\infty}\sum_{j=1}^n \binom nj {(-1)^{j-1}\over (j-1)!}
\sum_{m<X} {\Lambda(m)\over m^{3/2}} b_f(m)(\ln m)^{j-1}$$ 
which exists. \endproclaim

   \proclaim{Lemma 3.2} (see [14])
Let $F(x)$ be a function defined on $\Bbb R$ such that
$$2F(x)=F(x+0)+F(x-0)$$
for all $x\in\Bbb R$, such that $F(x)\exp((\epsilon+1/2)|x|)$
is integrable and of bounded variation on $\Bbb R$ for a
constant $\epsilon >0$, and such that
$(F(x)-F(0))/x$ is of bounded variation on $\Bbb R$.   Then
$$\aligned \sum_\rho\Phi(\rho)
= &2F(0)\ln{\sqrt N\over 2\pi}-\sum_{n=1}^\infty
{\Lambda(n)\over n}b_f(n)[F(\ln n)+F(-\ln n)]\\
&-\int_0^\infty \left({F(x)+F(-x)\over e^x-1}-2F(0)
{e^{-x}\over x}\right)dx,\endaligned $$
where the sum on $\rho$ runs over all zeros of $\xi_f(s)$ in
the order given by $|\Im\rho|<T$ for $T\to\infty$, and
$$\Phi(s)=\int_{-\infty}^\infty F(x) e^{(s-1/2)x}dx. $$
\endproclaim

\proclaim{Lemma 3.3} ([5] [6] [15])  Let $f$ be a newform of weight $2$
for $\Gamma_0(N)$.  Then there an absolute effective constant $c>0$
such that $L_f(s)$ has no zeros in the region
$$\{s=\sigma+it:\, \sigma\geq 1-{c\over \ln(N+1+|t|)}\},$$
where $L_f(s)$ is given in (1.2).
\endproclaim

 \proclaim{Lemma 3.4}  (Lemma 2 of [2])  For $n=1,2,\cdots$, let
$$F_n(x)=\cases e^{x/2}\sum_{j=1}^n \binom nj {x^{j-1}\over (j-1)!},
&\text{if $-\infty<x<0$};\\
n/2, &\text{if $x=0$};\\
0, &\text{if $0<x$}. \endcases $$
Then
$$\Phi_n(s)=1-\left(1-{1\over s}\right)^n$$
where $\Phi_n$ is related to $F_n$ by the relation
$$\Phi_n(s)=\int_{-\infty}^\infty F_n(x) e^{(s-1/2)x}dx. $$
\endproclaim

  \demo{Proof of Theorem 3.1}
 Since $\xi_f(s)$ is an entire function of order one
  and satisfies the functional identity
  $\xi_f(s)=w\xi_f(1-s)$,  we have
$$\xi_f(s)=w\xi_f(1)\prod_\rho (1- s/\rho)$$
where the product is over all the zeros of $\xi_f(s)$ in
the order given by $|\Im\rho|<T$ for $T\to\infty$.
If $\varphi_f(z)=\xi_f(1/(1-z))$, then
$${\varphi_f^\prime(z)\over\varphi_f(z)}=\sum_{n=0}^\infty
\tau_f(n+1)z^n \tag 3.3$$
where the coefficients $\tau_f(n)$ are given in (3.1).

    For a sufficiently large positive number $X$
that is not an integer, let
$$F_{n,X}(x)=\cases F_n(x), &\text{if $-\ln X<x<\infty$};\\
{1\over 2}F_n(-\ln X), &\text{if $x=-\ln X$}; \\
0, &\text{if $-\infty<x<-\ln X$} \endcases $$
where $F_n(x)$ is given in Lemma 3.4.  Then $F_{n,X}(x)$
satisfies all conditions of Lemma 3.2.  Let
$$\Phi_{n, X}(s)=\int_{-\infty}^\infty F_{n, X}(x)
e^{(s-1/2)x}dx.$$
By Lemma 3.2, we obtain that
$$\aligned \sum_\rho \Phi_{n,X}(\rho)
= &2F_{n,X}(0)\ln {\sqrt N\over 2\pi}
-\sum_{k=1}^\infty {\Lambda(k)\over k} b_f(k)
[F_{n,X}(\ln k)+F_{n,X}(-\ln k)]\\
&-\int_0^\infty \left({F_{n,X}(x)+F_{n,X}(-x)\over e^x-1}-2F_{n,X}(0)
{e^{-x}\over x}\right)dx,\endaligned $$
where the sum on $\rho$ runs over all zeros of $\xi_f(s)$ in the
order given by $|\Im\rho|<T$ for $T\to\infty$.
It follows that
$$\aligned \lim_{X\to\infty}&\sum_\rho\Phi_{n,X}(\rho)\\
&= n\left(\ln {\sqrt N\over 2\pi}-\gamma\right)
-\lim_{X\to\infty}\sum_{j=1}^n \binom nj {(-1)^{j-1}\over (j-1)!}
\sum_{k<X}{\Lambda(k)\over k^{3/2}} b_f(k)(\ln k)^{j-1}\\
&+n\left(-{2\over 3}+\sum_{l=1}^\infty {3\over l(2l+3)}\right)
+\sum_{j=2}^n\binom nj(-1)^j\sum_{l=1}^\infty
{1\over (l+1/2)^j}. \endaligned \tag 3.4 $$
 We have
$$\aligned \Phi_n(s)-\Phi_{n,X}(s)&=X^{-s}\sum_{j=1}^n\binom nj
(-1)^{j-1}\sum_{k=0}^{j-1}{(\ln X)^{j-k-1}\over (j-k-1)!}
 s^{-k-1}\\
 &={X^{-s}\over s}\sum_{j=1}^n\binom nj {(-\ln X)^{j-1}\over (j-1)!}
 +O\left({(\ln X)^{n-2}\over |s|^2}
 X^{-\Re s}\right). \endaligned  \tag 3.5$$

   Let $\rho$ be any zero of $\xi_f(s)$.
 By Lemma 3.3, we have
 $${c\over \ln (N+1+|\rho|) }
 \leq \Re\rho\leq 1-{c\over \ln (N+1+|\rho|)}  $$
 for a positive constant $c$.  An argument similar to that
 made in the proof of (3.9) of [2] shows that
 $$\sum_\rho {X^{-\Re\rho}\over |\rho|^2}
 \ll e^{-c^\prime\sqrt{\ln X}} \tag 3.6$$
 for a positive constant $c^\prime$.

    Since
     $$\aligned \sum_\rho {X^{-\rho}\over \rho}&=\sum_\rho
   {X^{-(1-\rho)}\over 1-\rho} \\
   &=-{1\over X}\sum_\rho {X^\rho\over\rho}
   +O\left(\sum_\rho {X^{-(1-\Re\rho)}\over |\rho|^2}\right)\\
   &=-{1\over X}\sum_\rho {X^\rho\over\rho}
 +O\left(e^{-c^\prime\sqrt{\ln X}}\right),
   \endaligned $$
and since
$$\lim_{X\to\infty}{(\ln X)^{j-1}\over X}\sum_\rho {X^\rho\over\rho}=0 $$
for $j=1,2,\cdots, n$ by Theorem 4.2 and Theorem 5.2 of [15],
we have
$$\lim_{X\to\infty}(\ln X)^{j-1}\sum_\rho {X^{-\rho}\over \rho}=0 \tag 3.7$$
 for $j=1,2,\cdots, n$.  It follows from (3.5), (3.6) and
 (3.7) that
  $$\lim_{X\to\infty}\sum_\rho\Phi_{n,X}(\rho)
  =\sum_\rho\Phi_n(\rho). \tag 3.8$$
Since $\xi_f(s)$ is an entire function of order one, the
series 
$$\sum_\rho\Phi_n(\rho)$$
is convergent, where the sum on $\rho$ runs over all zeros of $\xi_f(s)$ 
in the order given by $|\Im\rho|<T$ for $T\to\infty$.
Hence, by (3.4) and (3.8) the limit  
$$\lim_{X\to\infty}\sum_{j=1}^n \binom nj {(-1)^{j-1}\over (j-1)!}
\sum_{k<X}{\Lambda(k)\over k^{3/2}} b_f(k)(\ln k)^{j-1}$$
exists.

This completes the proof of the theorem.
\qed\enddemo

\heading
4.  Preliminary results
\endheading

  In this section, we collect some technical results for
the proof of Theorem 1.2.

  A fundamental result of Hecke asserts that a basis
$\{f_1, f_2, \cdots, f_g\}$ in $S_2(N)$ exists which consists
of eigenfunctions of all the Hecke operators $T(n)$ with
$(n, N)=1$;  see  Theorem 6.21 in Iwaniec [8].   We can assume
that each $f_j$ is either a normalized newform in $S_2(N)$ or
coming from a normalized newform in a lower level.
For $j=1,\cdots, g$,
we choose $g_j=f_j$ if $f_j$ is a normalized newform in $S_2(N)$,
and $g_j=f^\prime_j$ if $f_j$ is an old form in $S_2(N)$
and if $f^\prime_j$ is a normalized newform in $S_2(N^\prime)$ for
some divisor $N^\prime$ of $N$ such that $f_j(z)=f^\prime_j(dz)$
for some positive integer $d|N/N^\prime$.

   Let
$$\xi_H(s)=\prod_{j=1}^g \xi_{g_j}(s) \tag 4.1$$
where $\xi_{g_j}(s)$ is defined as in (1.3).
Since $\xi_{g_j}(s)$ is an entire function and satisfies
the functional identity $\xi_{g_j}(s)=w_j \xi_{g_j}(1-s)$,
where $w_j=(-1)^{r_j}$ with $r_j$ being the vanishing
order of $\xi_{g_j}$ at $s=1/2$, the function $\xi_H(s)$ is entire
and satisfies the functional identity
$$\xi_H(s)=\epsilon \xi_H(1-s), \tag 4.2$$
where $\epsilon =\pm 1$.  Put
$$\tau_H(n)=\sum_{j=1}^g\tau_{g_j}(n),\tag 4.3$$
where $\tau_{g_j}(n)$ is defined similarly as in (3.1).
If  $\varphi_H(z)=\xi_H(1/(1-z))$, then we have
$${\varphi_H^\prime(z)\over\varphi_H(z)}=\sum_{n=0}^\infty
\tau_H(n+1)z^n \tag 4.4$$
by (3.3).

  \proclaim{Lemma 4.1}  For all positive integers $n$ with
  $(n, N)=1$, we have
    $$(T(n)f_j)(z)=\lambda_{g_j}(n)f_j(z)$$
    for $j=1,2,\cdots, g$, where $\lambda_{g_j}(n)$ is
  the eigenvalue of $T(n)$ acting on $g_j(z)$. \endproclaim

\demo{Proof}  If $f_j$ is a newform, the stated identity is
trivially true.

   Next, we assume that $f_j$ is an old form.
Let $g_j$ be a normalized newform in $S_2(N^\prime)$
for some divisor $N^\prime$ of $N$ such that $f_j(z)=g_j(dz)$
for some positive integer $d|N/N^\prime$.  Since $(n,N)=1$,
by (1.1) we have
$$(T(n) f)(z)={1\over n}\sum_{ad=n}a^2
\sum_{0\leq b<d} f\left({az+b\over d}\right)$$
for any function $f$ in $S_2(N)$ or $S_2(N^\prime)$.
Thus, we have
$$(T(n) f_j)(z)={1\over n}\sum_{\alpha\delta=n}\alpha^2
\sum_{0\leq \beta<\delta} g_j
\left({\alpha dz+d\beta\over \delta}\right). $$
Since $(n, N)=1$, $d|N$ and $\alpha\delta=n$, we have
$(\delta, d)=1$.  Let $r_\beta$ be remainder of
$d\beta$ modulo $\delta$.  Then $\{r_\beta: 0\leq\beta<\delta\}
=\{0, 1, \cdots, \delta-1\}$.  Since
$g_j\in S_2(N^\prime)$, we have
$$g_j\left({\alpha dz+d\beta\over \delta}\right)
=g_j \left({\alpha dz+r_\beta\over \delta}\right).$$
It follows that
$$\aligned (T(n) f_j)(z)&={1\over n}\sum_{\alpha\delta=n}\alpha^2
\sum_{0\leq \beta<\delta} g_j
\left({\alpha dz+\beta\over \delta}\right)\\
&=(T(n) g_j)(w)=\lambda_{g_j}(n)g_j(w)
=\lambda_{g_j}(n)f_j(z) \endaligned  $$
where $w=dz$.

  This completes the proof of the lemma.
\qed\enddemo

  \proclaim{Lemma 4.2}  Let $\xi_N(s)$ be given in (1.5).
  Then $\xi_N(s)$ is an entire function, and its zeros
 in the strip $0<\Re s<1$ appear in pairs $\rho$ and $1-\rho$.
 \endproclaim

  \demo{Proof}    For $j=1,2, \cdots, g$, if $g_j$ is a normalized
  newform in $S_2(N_j)$, then we can write
  $$L_{g_j}(s)=\prod_{p|N_j}(1-\lambda_{g_j}(p)p^{-s})^{-1}
\prod_{p\nmid N_j} (1-\lambda_{g_j}(p)p^{-s}+p^{1-2s})^{-1}.
\tag 4.5$$
Since $f_1, f_2, \cdots, f_g$ are eigenfunctions of all the
Hecke operators $T(n)$ with $(n, N)=1$, by Lemma 4.1 we have
$$\det|1-T(p)p^{-s}+p^{1-2s}I|=\prod_{j=1}^g
(1-\lambda_{g_j}(p)p^{-s}+p^{1-2s})$$
for any prime $p\nmid N$.  Let $L_N(s)$ be given in (1.4).
Then we have
$$\aligned L_N(s)&=\prod_{j=1}^g\prod_{p\nmid N}
(1-\lambda_{g_j}(p)p^{-s}+p^{1-2s})^{-1}\\
&=\prod_{j=1}^g\left(L_{g_j}(s)\prod_{p|N_j}(1-\lambda_{g_j}(p)p^{-s})
\prod_{p\nmid N_j,p|N}(1-\lambda_{g_j}(p)p^{-s}+p^{1-2s})\right).
\endaligned $$
Since
$$\xi_{g_j}(s)=N_j^{s/2}(2\pi)^{-s}
\Gamma({1\over 2}+s)L_{g_j}({1\over 2}+s),$$
we have
$$\xi_H(s)=A^{s/2}N^{gs/2}(2\pi)^{-gs}\Gamma^g
\left({1\over 2}+s\right)\prod_{j=1}^gL_{g_j}
\left({1\over 2}+s\right)$$
where $A=N^{-g}\prod_{j=1}^g N_j$.  It follows that
$$\xi_N(s)=\xi_H(s) A^{-s/2}
\prod_{j=1}^g(\prod_{p|N_j}(1-\lambda_{g_j}(p)p^{-s-1/2})
\prod_{p\nmid N_j,p|N}(1-\lambda_{g_j}(p)p^{-s-1/2}+p^{-2s})).
\tag 4.6$$
This implies that $\xi_N(s)$ is an entire function.
Since $|\lambda_{g_j}(p)|=\sqrt p$ for $p|N_j$
by ii) in Theorem 3 of Li [10] and since
the two roots of the polynomial
$1-\lambda_{g_j}(p)p^{-1/2}z+z^2$ for $p\nmid N_j$
are conjugate complex numbers of absolute value one by
the Ramanujan conjecture which was proved in Th\'eor\`eme 8.2 of
Deligne [4], zeros of $\xi_N(s)$ in the critical strip $0<\Re s<1$
 appear in pairs $\rho, 1-\rho$ by (4.2) and (4.6).

  This completes the proof of the lemma.
\qed\enddemo

\demo{Proof of Corollary 1.4}  Let $\tau_N(n)$ be defined by
(1.6) for all positive integers $n$.
Since $\xi_H(s)$ is an entire function of order one,
by (4.6) $\xi_N(s)$ is an entire function of order one.
This implies that
$$\sum_\rho {1+|\Re \rho|\over (1+|\rho|)^2}<\infty,$$
where the sum is over all zeros $\rho$ of $\xi_N(s)$.
Thus, conditions of Theorem 1.3 are satisfied.
Since by (4.6) all zeros of $\xi_N(s)$ outside the
strip $0<\Re s<1$ lie on the line $\Re s=0$,
 Theorem 1.3 implies that all zeros of $\xi_N(s)$ in the
 critical strip $0<\Re s<1$ satisfy
$\Re s\leq 1/2$ if, and only if, $\tau_N(n)\geq 0$
for all positive integers $n$.  By Lemma 4.2, all zeros of
of $\xi_N(s)$ in the  critical strip $0<\Re s<1$
appear in pairs $\rho$ and $1-\rho$.
Thus, $\tau_N(n)\geq 0$
for all positive integers $n$ if, and only if,
$\Re \rho\leq 1/2$ and $\Re (1-\rho)\leq 1/2$
for all zeros $\rho$ of $\xi_N(s)$ in the
critical strip $0<\Re s<1$.  That is,
all zeros of $\xi_N(s)$ in the
 strip $0<\Re s<1$ lie on the critical line
$\Re s=1/2$ if, and only if, $\tau_N(n)\geq 0$
for all positive integers $n$.

 This completes the proof of the corollary.
\qed\enddemo

\proclaim{Lemma 4.3}  Let $p$ be a prime, and let
$\alpha$ be a complex number of absolute value one.
Then we have
$$1-\alpha p^{-s}=c_p s^{\epsilon_p}\prod_\rho (1-s/\rho)$$
where the product on $\rho$ is over all nonzero zeros
of $1-\alpha p^{-s}$
taken in the order given by $|\rho|<T$ for $T\to\infty$ and where
$c_p=1-\alpha, \epsilon_p =0$ if $\alpha\neq 1$ and
$c_p=\ln p, \epsilon_p=1$ if $\alpha=1$.
\endproclaim

 \demo{Proof}  Since $1-\alpha p^{-s}$ is an entire function of
 order one, by Hadamard's factorization theorem there
 is a constant $a$ such that
 $$1-\alpha p^{-s}=c_pe^{as}s^{\epsilon_p}
 \prod_\rho (1-s/\rho)e^{s/\rho} \tag 4.7$$
 where the product is over all nonzero zeros of $1-\alpha p^{-s}$.
 Let $\alpha=e^{it}$ with $0\leq t<2\pi$.
 Then the zeros of $1-\alpha p^{-s}$
are $i(t+2k\pi)/\ln p$, $k=0, \pm 1, \pm 2, \cdots$.
Since
$$\sum_{k=1}^\infty \left({-i\ln p\over t+2k\pi}
+{-i\ln p\over t-2k\pi}\right)
=2it\ln p\sum_{k=1}^\infty {1\over (2k\pi)^2-t^2} $$
is absolutely convergent,  by using (4.7) we can write
$$1-\alpha p^{-s}=c_pe^{hs}s^{\epsilon_p}
 \prod_\rho (1-s/\rho) \tag 4.8 $$
 for a constant $h$, where the product runs over all
 nonzero zeros $\rho$ of $1-\alpha p^{-s}$
taken in the order given by $|\rho|<T$ for $T\to\infty$.
By taking logarithmic derivative of both sides of (4.8)
with respect to $s$ we get
$${\alpha\ln p \over p^s-\alpha }
=h+{\epsilon_p \over s}+\sum_\rho {1\over s-\rho}.
\tag 4.9$$
By letting $s\to\infty$ in (4.9) we find that $h=0$.
Then the stated identity follows from (4.8).

  This completes the proof of the lemma.
\qed\enddemo

  \proclaim{Lemma 4.4}  Let $\alpha, p$ be given in
  Lemma 4.3, and let $s=(1-z)^{-1}$.  Then we have
  $$ {d\over dz}\ln (1-\alpha p^{-s})
  =\sum_{n=0}^\infty \left(\sum_\rho
  [1-(1-1/\rho)^{-n-1}]\right)z^n$$
for $z$ in a small neighborhood of the origin,
where the sum on $\rho$ is over all zeros
of $1-\alpha p^{-s}$
taken in the order given by $|\rho|<T$ for $T\to\infty$.
\endproclaim

 \demo{Proof} By Lemma 4.3 we have
$${d\over dz}\ln (1-\alpha p^{-s})={\epsilon_p\over 1-z}
+\sum_\rho {1\over (1-\rho)+\rho z}{1\over 1-z} \tag 4.10$$
where the sum on $\rho$ is over all nonzero zeros of
$1-\alpha p^{-s}$.  Since
$$\aligned {1\over (1-\rho)+\rho z}{1\over 1-z}
&={1\over 1-\rho}\{\sum_{k=0}^\infty \left({-\rho\over
1-\rho}\right)^kz^k\}\{\sum_{l=0}^\infty z^l\}\\
&=\sum_{n=0}^\infty [1-(1-1/\rho)^{-n-1}] z^n
\endaligned $$
for $z$ in a small neighborhood of the origin,
the stated identity follows from (4.10).
Note that, if $s=\rho=0$ is
a zero of $1-\alpha p^{-s}$, then $(1-1/\rho)^{-n}$ is
interpreted to be $0$ for all positive integers $n$.

  This completes the proof of the lemma.
\qed\enddemo

    Remark 4.5.  By (4.6), Lemma 4.4, (3.3) and (4.1) we have
$${d\over dz} \log \left(A^{s/2}\xi_N(s)\right)
=\sum_{n=0}^\infty \tau_N(n+1) z^n  $$
with $s=(1-z)^{-1}$, where the $\tau_N(n)$'s are given
in (1.6).   Then, by Lemma 4.2, to prove
that all zeros of $\xi_N(s)$ in the strip $0<\Re s<1$
lie on the critical line $\Re s=1/2$ it is enough to
find an upper bound for each $\tau_N(n)$ which implies
that the above series is analytic for $|z|<1$.

  \proclaim{Lemma 4.6}  Let $\alpha, p$ be given in
  Lemma 4.3, and let $s=(1-z)^{-1}$.  Then we have
  $$ {d\over dz}\ln (1-\alpha p^{-s})
  =\sum_{n=0}^\infty \left(\sum_{j=0}^n \binom {n+1}{j+1}
  {(-1)^j\over j!}\sum_{k=1}^\infty {\ln p\over p^k}\alpha^k
  (\ln p^k)^j \right)z^n$$
for $z$ in a small neighborhood of the origin.
\endproclaim

 \demo{Proof}  Let $m$ be any positive integer.  By using
 mathematical induction on $n$, we can show that
 $${d^n\over dz^n}\left[(1-z)^{-m}p^{-ks}\right]_{|z=0}
 =p^{-k}\sum_{j=0}^n\binom nj (n+m-1)\cdots (j+m)
 (-\ln p^k)^j \tag 4.11 $$
 for $n=1,2,\cdots$.  By using the formula (4.11) with
 $m=2$ we find that
 $$\aligned {d\over dz}\ln (1-\alpha p^{-s})
  &=(1-z)^{-2}\ln p\sum_{k=1}^\infty \alpha^k p^{-ks}\\
  &=\sum_{n=0}^\infty {z^n\over n!}\left(\sum_{k=1}^\infty
 \alpha^k {\ln p\over p^k}\sum_{j=0}^n\binom nj
 (n+1)\cdots (j+2)(-\ln p^k)^j \right)\\
 &=\sum_{n=0}^\infty z^n\left(\sum_{j=0}^n \binom {n+1}{j+1}
  {(-1)^j\over j!}\sum_{k=1}^\infty {\ln p\over p^k}\alpha^k
  (\ln p^k)^j \right).
 \endaligned $$

  This completes the proof of the lemma.
\qed\enddemo

\proclaim{Lemma 4.7}  Let $p$ be a prime, and let
$\alpha$ be a complex number of absolute value one.
Then we have
$$\sum_\rho [1-(1-1/\rho)^{-n}]
  =\sum_{j=1}^n \binom nj {(-1)^{j-1}\over (j-1)!}
\sum_{k=1}^\infty {\ln p\over p^k}\alpha^k(\ln p^k)^{j-1}
  $$
for $n=1,2,\cdots$, where the sum on $\rho$ is over
all zeros of $1-\alpha p^{-s}$ taken in the order given by
$|\Im \rho|<T$ for $T\to\infty$ with a zero of multiplicity
$\ell$ appearing $\ell$ times in the list.
\endproclaim

\demo{Proof}  The stated identity follows from
Lemma 4.4 and Lemma 4.6. \qed\enddemo

  \proclaim{Lemma 4.8}  For $n=1,2,\cdots$, we have
$$\tau_N(n)=\sum_{j=1}^g\tau_{g_j}(n)
+\sum_{l=1}^n \binom nl {(-1)^{l-1}\over (l-1)!}
\sum_{(m, N)>1, m=1}^\infty {\Lambda(m)\over m^{3/2}}
(\sum_{j=1}^g b_{g_j}(m))(\ln m)^{l-1} $$
where $b_{g_j}(m)$ is given as in (3.2).
\endproclaim

  \demo{Proof}  Let $\tau_H(n)$ be given in (4.3).
By (3.1), (4.1) and (4.6) we have
  $$\tau_N(n)=\tau_H(n)+\sum_{j=1}^g \{
  \sum_{p|N_j}\sum_{\rho_j}[1-(1-1/\rho_j)^{-n}]
  +\sum_{p\nmid N_j, p|N}\sum_{\alpha_j}
  [1-(1-1/\alpha_j)^{-n}]\},\tag 4.12 $$
  where the sum on $\rho_j$ is over all zeros of
$1-\lambda_{g_j}p^{-s-1/2}$ with $p|N_j$ and where the sum
on $\alpha_j$ is over all zeros of
$1-\lambda_{g_j}(p)p^{-s-1/2}+p^{-2s}$
 with $p\nmid N_j, p|N$.  By Lemma 4.7 we have
 $$\sum_{\rho_j}[1-(1-1/\rho_j)^{-n}]
 =\sum_{l=1}^n \binom nl {(-1)^{l-1}\over (l-1)!}
\sum_{k=1}^\infty {\ln p\over p^{3k/2}}
b_{g_j}(p^k)(\ln p^k)^{l-1} \tag 4.13 $$
and
$$\sum_{\alpha_j}[1-(1-1/\alpha_j)^{-n}]
=\sum_{l=1}^n \binom nl {(-1)^{l-1}\over (l-1)!}
\sum_{k=1}^\infty {\ln p\over p^{3k/2}}
b_{g_j}(p^k)(\ln p^k)^{l-1}. \tag 4.14 $$
It follows from (4.13) and (4.14) that
$$\aligned &\sum_{p|N_j}\sum_{\rho_j}
  [1-(1-1/\rho_j)^{-n}]+\sum_{p\nmid N_j, p|N}\sum_{\alpha_j}
  [1-(1-1/\alpha_j)^{-n}]\\
  =&\sum_{p|N}\sum_{l=1}^n \binom nl {(-1)^{l-1}\over (l-1)!}
\sum_{k=1}^\infty {\ln p\over p^{3k/2}}
b_{g_j}(p^k)(\ln p^k)^{l-1}. \endaligned \tag 4.15 $$
By (4.15) we have
$$\aligned &\sum_{j=1}^g\left(\sum_{p|N_j}\sum_{\rho_j}
  [1-(1-1/\rho_j)^{-n}]+\sum_{p\nmid N_j, p|N}\sum_{\alpha_j}
  [1-(1-1/\alpha_j)^{-n}]\right)\\
  &=\sum_{l=1}^n \binom nl {(-1)^{l-1}\over (l-1)!}
\sum_{(m, N)>1, m=1}^\infty {\Lambda(m)\over m^{3/2}}
(\sum_{j=1}^g b_{g_j}(m))(\ln m)^{l-1}.\endaligned $$
The stated identity then follows from (4.12).

  This completes the proof of the lemma.
\qed\enddemo

\heading
5.  Proof of Theorem 1.2
\endheading

  In this section we complete the proof of Theorem 1.2.

   We define an operator $S$ acting on the space $S_2(N)$ by
  $$\aligned &S(1)=2I \\
  & S(p)=T(p) \\
  &S(p^m)=T(p^m)-pT(p^{m-2})\endaligned  \tag 5.1$$
for $m=2,3,\cdots$; see Ihara [7].

  \proclaim{Lemma 5.1}  For each prime $p\nmid N$, the trace
$\text{tr}(S(p^m))$ of $S(p^m)$ acting on the space $S_2(N)$
is given by
$$\text{tr}(S(p^m))=\sum_{j=1}^g b_{g_j}(p^m) \tag 5.2 $$
for $m=0,1,2,\cdots$.
\endproclaim

\demo{Proof}  Let $p$ be any prime.
 It follows from the recursion formula (see (6.25) in Iwaniec [8])
$$T(p^{m+1})=T(p)T(p^m)-pT(p^{m-1})$$
that
$$S(p^{m+1})=S(p)S(p^m)-pS(p^{m-1}) \tag 5.3 $$
for $m=1,2,\cdots$.

  Let $p$ is any prime not dividing $N$.  When $m=0$, we have
$$S(p^m)f_j=2f_j=b_{g_j}(p^m)f_j.$$
When $m=1$, we have
$$S(p^m)f_j=T(p)f_j=\lambda_{g_j}(p)f_j=b_{g_j}(p^m)f_j$$
by Lemma 4.1.  Assume that
$$S(p^m)f_j=b_{g_j}(p^m)f_j$$
for all integers $m\leq k$.  When $m=k+1$, we have
$$\aligned S(p^m)f_j&=(S(p)S(p^k)-pS(p^{k-1}))f_j\\
&=(b_{g_j}(p)b_{g_j}(p^k)-pb_{g_j}(p^{k-1}))f_j\\
&=b_{g_j}(p^m)f_j.\endaligned $$
By mathematical induction the identity
$$S(p^m)f_j=b_{g_j}(p^m)f_j$$
holds for all nonnegative integers $m$.
Since $\{f_1, \cdots, f_g\}$ is a basis for $S_2(N)$, we have
$$tr(S(p^m))=\sum_{j=1}^g b_{g_j}(p^m)$$
for $m=0, 1,2,\cdots$.

  This completes the proof of the lemma.
\qed\enddemo

  From Lemma 5.1 and the definition (5.1) we
obtain the following corollary.

   \proclaim{Corollary 5.2}  For each prime $p\nmid N$, we have
$$\sum_{j=1}^g b_{g_j}(p^m)=\text{tr}(T(p^m))-p\,\text{tr}(T(p^{m-2})$$
for $m=0,1,2,\cdots$.
\endproclaim

 \demo{Proof of Theorem 1.2}  By Lemma 4.8 we have
$$\tau_N(n)=\sum_{j=1}^g\tau_{g_j}(n)
+\sum_{l=1}^n \binom nl {(-1)^{l-1}\over (l-1)!}
\sum_{(m, N)>1, m=1}^\infty {\Lambda(m)\over m^{3/2}}
(\sum_{j=1}^g b_{g_j}(m))(\ln m)^{l-1} \tag 5.4 $$
for $n=1,2,\cdots$.  Since $g_j$ is a normalized
newform in $S_2(N_j)$, by Theorem 3.1 we have
$$ \aligned \tau_{g_j}(n)=& n\left(\ln {\sqrt {N_j}\over 2\pi}
-\gamma\right)
-\sum_{l=1}^n \binom nl {(-1)^{l-1}\over (l-1)!}\sum_{m=1}^\infty
{\Lambda(m)\over m^{3/2}} b_{g_j}(m)(\ln m)^{l-1}\\
&+n\left(-{2\over 3}+\sum_{l=1}^\infty {3\over l(2l+3)}\right)
+\sum_{m=2}^n\binom nm (-1)^m\sum_{l=1}^\infty
{1\over (l+1/2)^m} \endaligned \tag 5.5  $$
By using (5.4) and (5.5) we obtain that
$$\aligned \tau_N(n)=& {n\over 2}\ln(N_1\cdots N_g)
-\sum_{l=1}^n \binom nl {(-1)^{l-1}\over (l-1)!}
\sum_{\underset{(m, N)=1}\to {m=1}}^\infty {\Lambda(m)\over m^{3/2}}
\{\sum_{j=1}^g b_{g_j}(m)\}(\ln m)^{l-1}\\
&-ng\left(\ln 2\pi +\gamma+{2\over 3}
-\sum_{l=1}^\infty {3\over l(2l+3)}\right)
+g\sum_{m=2}^n\binom nm \sum_{l=1}^\infty
{(-1)^m\over (l+1/2)^m}. \endaligned  \tag 5.6$$
By Theorem 5 of Atkin and Lehner [1], we have
$$\prod_{N_j\neq N, 1\leq j\leq g}N_j
=\prod_{1<m<N, m|N} m^{\nu_m d(N/m)}.\tag 5.7$$
Since the dimension of the space $S_2(N)$ is $g$,
which is explicitly given in Proposition 1.40 and
Proposition 1.43 of Shimura [17],  we have the following
 recurrence formula for the number of newforms in a basis
 for $S_2(N)$:
$$\nu_N=g-\sum_{1<m<N, m|N}\nu_md(N/m). \tag 5.8 $$
Thus, we have
$$\prod_{N_j= N, 1\leq j\leq g}N_j=N^{\nu_N}. \tag 5.9$$
By (5.6), (5.7) and (5.9) we have
$$\aligned \tau_N(n)=& {n\over 2}
\ln\left(N^{\nu_N}\prod_{1<m<N, m|N} m^{\nu_md(N/m)}\right)\\
&-\sum_{l=1}^n \binom nl {(-1)^{l-1}\over (l-1)!}
\sum_{\underset{(m, N)=1}\to {m=1}}^\infty {\Lambda(m)\over m^{3/2}}
\{\sum_{j=1}^g b_{g_j}(m)\}(\ln m)^{l-1}\\
&-ng\left(\ln 2\pi +\gamma+{2\over 3}
-\sum_{l=1}^\infty {3\over l(2l+3)}\right)
+g\sum_{m=2}^n\binom nm \sum_{l=1}^\infty
{(-1)^m\over (l+1/2)^m}. \endaligned  \tag 5.10 $$
 Since
$$\ln 2=\sum_{n=1}^\infty {(-1)^{n-1}\over n},$$
 we have
$$\ln 2={4\over 3}-{3\over 2}\sum_{l=1}^\infty {1\over l(2l+3)}.$$
It follows that
$$\ln 2\pi +\gamma+{2\over 3}
-\sum_{l=1}^\infty {3\over l(2l+3)}=\ln(8\pi)+\gamma-2. \tag 5.11$$
Note that the identity (5.11) is due to B. Conrey.
The author wants to thank him for this observation.
By (5.10) and (5.11) we have
$$\aligned \tau_N(n)=& {n\over 2}
\ln\left(N^{\nu_N}\prod_{1<m<N, m|N} m^{\nu_md(N/m)}\right)\\
&-\sum_{l=1}^n \binom nl {(-1)^{l-1}\over (l-1)!}
\sum_{\underset{(m, N)=1}\to {m=1}}^\infty {\Lambda(m)\over m^{3/2}}
\{\sum_{j=1}^g b_{g_j}(m)\}(\ln m)^{l-1}\\
&-ng\left(\ln(8\pi)+\gamma-2\right)
+g\sum_{m=2}^n\binom nm \sum_{l=1}^\infty
{(-1)^m\over (l+1/2)^m}. \endaligned $$
The stated identity then follows from Corollary 5.2.

This completes the proof of the theorem.
\qed\enddemo

   The following result for the size of the last term 
on the right side of the identity in Theorem 1.2 is due to B. Conrey.
The author wants to thank him for allowing him
to include it here.

\proclaim{Lemma 5.3} (B. Conrey)  We have
$$\sum_{m=2}^n\binom nm \sum_{l=1}^\infty
{(-1)^m\over (l+1/2)^m}=n\ln n+O(n)$$
for all positive integers $n$.
\endproclaim

\demo{Proof}  We have
$$\aligned \sum_{m=2}^n\binom nm \sum_{l=1}^n{(-1)^m\over (l+1/2)^m}
&=\sum_{l=1}^n[\left(1-{1\over l+1/2}\right)^n-1+{n\over l+1/2}]\\
&=O(n)+n\sum_{l=1}^n{1\over l+1/2}\\
&=O(n)+n\ln n. \endaligned \tag 5.12 $$
Since
$$\sum_{l=n+1}^\infty {1\over (l+1/2)^m}<\int_n^\infty {1\over (t+1/2)^m}dt
={1\over m-1}(n+1/2)^{1-m}< n^{1-m},$$
we have
$$\left|\sum_{m=2}^n\binom nm
\sum_{l=n+1}^\infty{(-1)^m\over (l+1/2)^m}\right|
\leq  n(1+1/n)^n = O(n). \tag 5.13$$
The stated identity then follows from (5.12) and (5.13).

 This completes the proof of the lemma. \qed\enddemo

Remark 5.4.  It follows from Lemma 5.3 that the last term on 
right side of the identity in Theorem 1.2 is asymptotically
equal to $gn\ln n+O(n)$ as $n\to\infty$.
According to a recent result of Lagarias [9], we would have
$$\lim_{n\to\infty}{\tau_N(n)\over gn\ln n}=1$$
if all zeros of $\xi_N(s)$ in the strip $0<\Re s<1$ lie
on the critical line $\Re s=1/2$.

Remark 5.5.  The function $L_N(s)$ in (1.4) is a partial $L$-function
defined only for $p\nmid N$.  A question of the referee, which
author does not know, is how to define Euler factors for $p|N$ so
that the function $\xi_N(s)$ corresponding to a completed
$L$-function $L_N(s)$ satisfies a functional equation of automorphic
$L$-functions.

\enddocument